\documentstyle[12pt, euscript, fullpage]{amsart}
\renewcommand{\cal}{\EuScript}
\begin{document}
\newcommand{\scroll}[1]{\scrollmode#1\errorstopmode}
\newcommand{\comment}[1]{}

\newcommand{\createtitle}[2]{\title{#1}\author{Greg Martin}\address{Department of Mathematics\\University of Toronto\\Canada M5S 3G3}\email{gerg@@math.toronto.edu}\subjclass{#2}\maketitle}
\newcommand{\make}[1]{\label{#1sec}\noindent\input{#1}}
\newcommand{\Mmake}[1]{\label{#1sec}\noindent}


\newtheorem{theorem}{Theorem}
\newtheorem{lemma}[theorem]{Lemma}
\newtheorem{corollary}{Corollary}[theorem]
\newtheorem{proposition}[theorem]{Proposition}

\newenvironment{pflike}[1]{\noindent{\bf #1}}{\vskip10pt} 
\newenvironment{proof}{\begin{pflike}{Proof:}}{\qed\end{pflike}}

\def\(#1\)#2{
 \def\helpera{\ifcase#2(\or\big(\or\Big(\or\bigg(\or\Bigg(\else(\fi}
 \def\helperb{\ifcase#2)\or\big)\or\Big)\or\bigg)\or\Bigg)\else)\fi}
 \helpera{#1}\helperb}

\newcommand{\2}[1]{\ifmmode{\cal#1}\else$\cal#1$\fi}
\newcommand{\3}[1]{\#\{#1\}}
\newcommand{\abs}[1]{\left|#1\right|} 
\newcommand{\floor}[1]{\lfloor#1\rfloor}
\newcommand{\bfloor}[1]{\big\lfloor#1\big\rfloor}
\newcommand{\bbfloor}[1]{\bigg\lfloor#1\bigg\rfloor}
\newcommand{\ceil}[1]{\lceil#1\rceil}
\newcommand{\bceil}[1]{\big\lceil#1\big\rceil}
\newcommand{\bbceil}[1]{\bigg\lceil#1\bigg\rceil}

\renewcommand{\mod}[1]{{\ifmmode\text{\rm\ (mod~$#1$)}
 \else\discretionary{}{}{\hbox{ }}\rm(mod~$#1$)\fi}}
\newcommand{\ep}{\varepsilon}
\renewcommand{\implies}{\Rightarrow}
\newcommand{\rmif}{{\rm if\ }}

\newcommand{\half}{{\mathchoice{\textstyle\frac12}{1/2}{1/2}{1/2}}}
\newsymbol\dnd 232D 
\newcommand{\exdiv}{\mathrel{\mid\mid}}

\renewcommand{\lg}[1]{\mathop{\log_{#1}}}
\def\lgs#1^#2{\mathop{\log_{#1}^{#2}}}
\newcommand{\li}{\mathop{\rm li}}

\newcommand{\doublespace}{
  \baselineskip=24pt}
\newcommand{\spaceandahalf}{\parskip6pt\baselineskip=18pt}

\vfuzz=2pt 
\createtitle{Denser Egyptian Fractions}{11D68}

\newcommand{\lcm}{\mathop{\rm lcm}}

\section{Introduction}\label{introsec}\noindent An Egyptian fraction is a sum of reciprocals of distinct positive
integers, so called because the ancient Egyptians represented rational
numbers in that way. In an earlier paper, the author \cite{Mar:DEF}
showed that every positive rational number $r$ has Egyptian fraction
representations where the number of terms is of the same order of
magnitude as the largest denominator. More precisely, there exists a
positive constant $C(r)$ such that for every $x$ that is sufficiently
large in terms of $r$, there is a set \2E of positive integers not
exceeding $x$ with $\sum_{n\in\2E}1/n=r$ and $|\2E| > C(r)x$, where
$|\2E|$ denotes the cardinality of $\2E$.

Independently, Croot \cite{Cro:OSQoEaGaEF} considered the problem of
determining which integers can be represented as Egyptian fractions
with denominators not exceeding some bound $x$. He showed that any
integer less than
\begin{equation}
\log x + \gamma - \big( \frac92 + o(1) \big) {(\log\log x)^2\over\log x}
\label{nofx}
\end{equation}
(where $\gamma$ is Euler's constant) can be so represented, an
impressive feat considering that the sum of the reciprocals of all the
integers up to $x$ is only $\log x+\gamma$. The methods used in
\cite{Cro:OSQoEaGaEF} and \cite{Mar:DEF} are to an extent quite
coincident; one feature that they share is that the representations
constructed in both cases use only integers satisfying some upper
bound on the size of their prime factors (a restriction that at some
point is inherent in the problem rather than just a flaw in the
method). However, Croot's technical work in this regard is much
sharper than that in \cite{Mar:DEF}, as is attested to by the fact
that the expression (\ref{nofx}) is best possible (apart perhaps from
the constant $9/2$).

The purpose of this paper is to combine Croot's techniques with the
methods in the author's paper \cite{Mar:DEF} to sharpen the theorem
established in that paper, and to apply these methods to two related
questions concerning Egyptian fractions. We can improve the
aforementioned result from \cite{Mar:DEF} as follows:

\begin{theorem}
Let $r$ be a positive rational number. For every $x$ that is
sufficiently large in terms of $r$, there is a set\/ \2E of integers
not exceeding $x$, such that $\sum_{n\in\2E}1/n=r$ and
\begin{equation}
|\2E| > (1-e^{-r})x - O_r\big( {x\log\log x\over\log x} \big).
\label{Esize}
\end{equation}
Furthermore, this is best possible: the main term cannot be increased,
nor can the error term be reduced.
\label{denserthm}
\end{theorem}

Theorem \ref{denserthm} asserts the existence of Egyptian fraction
representations with many terms relative to a predetermined bound on
the maximal denominator; it is a slightly more delicate matter to turn
the problem around, and seek Egyptian fraction representations with
small denominators relative to a predetermined number of terms. To be
more concrete, for each positive integer $t$ let us define
\begin{equation*}
\2H_t(r) = \{ (x_1,\dots,x_t)\in{\bf Z}^t\colon x_1>\dots>x_t\ge1,\,
\sum_{i=1}^t \frac1{x_i}=r\},
\end{equation*}
the collection of all sets of denominators in $t$-term Egyptian
fraction representations of $r$, and
\begin{equation}
M_t(r) = \inf \{ x_1\colon (x_1,\dots,x_t)\in\2H_t(r) \},  \label{mrtdef}
\end{equation}
the smallest integer $x$ that is the largest denominator in a $t$-term
Egyptian fraction representation of $r$ (unless no such $t$-term
representation exists in which case $M_t(r)$ equals infinity). Among
the open problems posed by Erd\H os and Graham
\cite{ErdGra:OaNPaRiCNT} and also found in Guy \cite[Section
D11]{Guy:UPINT} is to determine the asymptotic size of $M_t(1)$, and
the question generalizes directly to~$M_t(r)$.

The set $\2H_t(r)$ might well be empty for small values of $t$. If we
define $t_0(r)$ to be the least number of terms in any Egyptian
fraction representation of $r$, then $\2H_t(r)$ is nonempty for every
$t\ge t_0(r)$, since a representation with $t$ terms can be converted
into one with $t+1$ terms by ``splitting'' the term with largest
denominator, using the identity
\begin{equation}
\frac1n = \frac1{n+1}+\frac1{n(n+1)}.  \label{split}
\end{equation}
Therefore $M_t(r)=\infty$ for $t<t_0(r)$ and $M_t(r)<\infty$ for $t\ge
t_0(r)$. (The only flaw in this argument arises when the largest
denominator is $n=1$, in which case the splitting identity does not
yield distinct terms; and in fact, the rational number 1 has an
Egyptian fraction representation with one term but no representation
with two terms. For this reason, we make the convention that
$t_0(1)=3$.)

The author's result from \cite{Mar:DEF} implies that for each positive
rational $r$, there are infinitely many values of $t$ for which
$M_t(r) \ll_r t$, and Theorem \ref{denserthm} improves the constant in
this upper bound to the best-possible value of
$(1-e^{-r})^{-1}+o(1)$. With a little more work, however, we need not
be content with establishing such a bound merely for infinitely
many~$t$. The following theorem completely resolves Erd\H os and
Graham's question in terms of the asymptotic behavior of $M_t(r)$:

\begin{theorem}
For all positive rational numbers $r$ and all integers $t\ge t_0(r)$,
we have
\begin{equation}
M_t(r) = \frac t{1-e^{-r}} + O_r\big( {t\log\log3t\over\log3t} \big),
\label{Msize}
\end{equation}
which is best possible.
\label{mrtthm}
\end{theorem}

Here it is appropriate to note that recently, Croot \cite{Cro:UFDSI}
established the related result that any positive rational $r$ has an
Egyptian fraction representation whose denominators all lie in the
interval $[x,e^rx+O_r(x\log\log x/\log x)]$, if $x$ is large enough in
terms of $r$. This result implies Theorem \ref{denserthm}, and in fact
somewhat more: if we define
\begin{equation*}
M'_t(r) = \min \{ x_t-x_1\colon (x_1,\dots,x_t)\in\2H_t(r) \},
\end{equation*}
then just as Theorem \ref{denserthm} implies that $M_t(r) =
t/(1-e^{-r}) + O_r(t\log\log t/\log t)$ for infinitely many $t$, so
does Croot's result imply the stronger result that $M'_t(r) = t +
O_r(t\log\log t/\log t)$ for infinitely many $t$. This asymptotic
expression for $M'_t(r)$ is also best possible; however, no result
analogous to Theorem \ref{mrtthm}, valid for all $t$, is obtained for
$M'_t(r)$.

The methods used to establish Theorems \ref{denserthm} and
\ref{mrtthm} allow us also to address another problem posed in
\cite{ErdGra:OaNPaRiCNT}. Erd\H os and Graham observe that a prime
power can never be the largest denominator in an Egyptian fraction
representation of 1, nor can a tiny multiple of a prime power%
; they ask whether the set of integers with this property has positive
density or even density 1. The analogous question can be asked about
those integers that cannot be the largest or second-largest
denominator in an Egyptian fraction representation of 1, and so on (as
is mentioned in \cite{Guy:UPINT}).

We can generalize this problem to Egyptian fraction representations of
any positive rational number $r$. For any positive integer $j$, let us
define
\begin{multline*}
\2L_j(r) = \{ x\in{\bf Z},\, x>r^{-1}\colon \hbox{there do not exist }
x_1,\dots,x_t\in{\bf Z},\, x_1>\dots>x_t\ge1 \\
\hbox{ with } { \textstyle\sum_{i=1}^t 1/x_i=r \hbox{ and }
x_j=x } \},
\end{multline*}
the set of numbers that {\it cannot\/} be the $j$th-largest
denominator in an Egyptian fraction representation of $r$. We exclude
the integers $x\le r^{-1}$ from consideration because they can never
be a denominator in an Egyptian fraction representation of $r$ (except
for the trivial representation when $r$ is itself the reciprocal of an
integer). The questions of Erd\H os and Graham then become whether
$\2L_1(1)$ has positive density, what can be said about
$\2L_1(1)\cap\2L_2(1)$, and so on.

In our primary theorem concerning these questions, we discover some
information about the sets $\2L_j(r)$ for $j\ge2$ that is perhaps
quite surprising:

\begin{theorem}
Let $r$ be a positive rational number. The set $\2L_j(r)$ is finite
for any integer $j\ge2$, and there exists an integer $j_0(r)$ such
that $\2L_j(r)$ is empty for all $j\ge j_0(r)$.
\label{finitethm}
\end{theorem}

One consequence of Theorem~\ref{finitethm} is that only finitely many
integers cannot be the second-largest denominator in an Egyptian
fraction representation of 1; possibly $\{2,4\}$ is a complete list of
integers (greater than 1) with this property. Another consequence is
that every integer greater than 1 can be the $j$th-largest denominator
in an Egyptian fraction representation of 1, when $j$ is sufficiently
large; possibly this holds for every $j\ge3$. It might be interesting
to establish a bound for the largest such integer and subsequently
determine the sets $\2L_2(1)$ and $\2L_3(1)$ precisely with the aid of
a computer; however, we do not undertake these tasks herein.

Because the prime factors of the denominator of $r$ are the only
primes that can possibly be the largest denominator in an
Egyptian fraction representation of $r$, the set $\2L_1(r)$ is
certainly infinite. However, we are able to answer Erd\H os and
Graham's question of whether $\2L_1(r)$ has positive density in the
negative; in fact, we can even establish the order of growth of
$\2L_1(r)$. Let us define the counting function $L_1(r;x)$ of
$\2L_1(r)$,
\begin{equation*}
L_1(r;x) = \3{1\le n\le x\colon n\in\2L_1(r)}.
\end{equation*}
Then we have the following theorem:

\begin{theorem}
Let $r$ be a positive rational number. The set $\2L_1(r)$ has zero
density, and in fact, if $x\ge3$ is a real number then
\begin{equation}
{x\log\log x\over\log x} \ll_r L_1(r;x) \ll_r {x\log\log x\over\log
x}.  \label{twosides}
\end{equation}
\label{lastthm}
\end{theorem}

The lower bound in the inequality (\ref{twosides}) is a simple
quantitative consequence of the observation that tiny multiples of
prime powers are elements of $\2L_1(r)$, as we shall show in the next
section, while the upper bound reflects the discovery that all
elements of $\2L_1(r)$ are of this form (the only ambiguity being the
exact meaning of ``tiny''). This discovery was probably known to Croot
(certainly in the case where $r$ is an integer), since the methods we
will use to establish the inequality (\ref{twosides}) are to a large
extent present in \cite{Cro:OSQoEaGaEF}.

The author would like to thank Ernest S.~Croot III for enlightening
conversations and for providing access to his manuscripts prior to
publication, as well as John Friedlander for suggestions that greatly
improved the presentation of this paper. The author also acknowledges
the support of National Science Foundation grant number DMS 9304580
and Natural Sciences and Engineering Research Council grant number
A5123.

\section{Reduction of Theorems~\ref{denserthm}
and~\ref{mrtthm}}\label{propssec}\noindent We begin by defining some notation that will be used throughout the
paper. Hereafter $p$ will always denote a prime and $q$ will always
denote a (not necessarily proper) prime power. As is standard, the
function $\pi(x)$ denotes the number of primes not exceeding $x$, and
$P(n)$ denotes the largest prime divisor of $n$. It is more convenient
for our purposes, however, to consider the sequence of prime powers as
more fundamental than the sequence of primes.  Therefore, we will use
$\pi^*(x)$ to denote the number of prime powers not exceeding $x$, so
that $\pi^*(x) = \pi(x) + \pi(x^{1/2}) + \pi(x^{1/3}) + \dotsb$. We
will also let $P^*(n)$ denote the largest prime power that divides
$n$; for example, $P(12)=3$ but $P^*(12)=4$.  (By convention we set
$P^*(1)=1$.) Notice that if $n$ is the least common multiple of $l$
and $m$, then $P^*(n)=\max\{P^*(l),P^*(m)\}$; in particular, if $a$
and $b$ are coprime and $a/b=a_1/b_1+a_2/b_2$, then
$P^*(b)\le\max\{P^*(b_1),P^*(b_2)\}$.

The statements of Theorems~\ref{denserthm} and~\ref{mrtthm} can be
reduced to the following two propositions:

\begin{proposition}
Let $I$ be a closed subinterval of $(0,\infty)$. There exists a
positive real number $T(I)$ such that, for all integers $t>T(I)$ and
all rational numbers $r=a/b\in I$ such that $P^*(b)<t\log^{-22}t$,
there is a set $\2E$ of $t$ distinct positive integers such that
$\sum_{n\in\2E}1/n=r$ and
\begin{equation}
\max\{n\in\2E\} < {t\over1-e^{-r}} + O_I\big( {t\log\log t\over\log t}
\big).
\label{OIsize}
\end{equation}
\label{estprop}
\end{proposition}

\begin{proposition}
Let $r$ be a positive rational number. There exists a positive
constant $\delta(r)$ such that, for every real number $x$ that is
sufficiently large in terms of $r$, all sets $\2E$ of positive
integers not exceeding $x$ for which $\sum_{n\in\2E}1/n=r$ satisfy
\begin{equation*}
|\2E| \le (1-e^{-r})x - \delta(r){x\log\log x\over\log x}.
\end{equation*}
\label{bestprop}
\end{proposition}

Let us first verify that Theorem~\ref{denserthm} follows from these
two propositions. Let $r=a/b$ be a positive rational number and $x$ a
real number that is sufficiently large in terms of $r$. Let $C(r)$ be
a large positive constant depending on $r$, and define $I$ to be the
length-zero interval $I=\{r\}$ and
\begin{equation}
t = \big\lceil (1-e^{-r})x - {C(r)x\log\log x\over\log x} \big\rceil.
\label{tdef}
\end{equation}
The integer $t$ will satisfy both $t>T(I)$ and $t\log^{-22}t>P^*(b)$
as soon as $x$ is large enough in terms of $r$ (since $T(I)$ depends
only on $r$), and so we can apply Proposition~\ref{estprop} to obtain
a set $\2E$ of $t$ distinct positive integers such that
$\sum_{n\in\2E}1/n=r$ and
\begin{equation*}
\max\{n\in\2E\} < {t\over1-e^{-r}} + O_r\big( {t\log\log t\over\log t}
\big).
\end{equation*}
But if $C(r)$ is chosen large enough, the right-hand side of this
inequality will be bounded above by $x$ by the definition~(\ref{tdef})
of $t$, and so the integers in \2E do not exceed $x$. Also, since
$|\2E|=t$, the definition of $t$ certainly implies the lower
bound~(\ref{Esize}), and therefore the first assertion of
Theorem~\ref{denserthm} has been shown to follow from
Proposition~\ref{estprop}. On the other hand, the remaining assertion
that the lower bound~(\ref{Esize}) is best possible follows
immediately from Proposition~\ref{bestprop}.

In the same manner we can deduce Theorem~\ref{mrtthm} from
Propositions~\ref{estprop} and~\ref{bestprop}. Given a positive
rational number $r=a/b$ and an integer $t>T(\{r\})$ that is so large
that $t\log^{-22}t>P^*(b)$, we can apply Proposition~\ref{estprop} to
obtain a set \2E of $t$ distinct positive integers such that
$\sum_{n\in\2E}1/n=r$, whose largest term is at most $t/(1-e^{-r}) +
O_r(t\log\log t/\log t)$. By the definition~(\ref{mrtdef}) of
$M_t(r)$, this establishes the upper bound implicit in the asymptotic
formula~(\ref{Msize}) when $t$ is sufficiently large in terms of $r$;
but by adjusting the constant implicit in the $O$-notation if
necessary, we see that this upper bound is valid for all $t\ge
t_0(r)$. On the other hand, if \2E is any set of $t$ distinct positive
integers satisfying $\sum_{n\in\2E}1/n=r$ whose largest element is
$x_1$, then Proposition~\ref{bestprop} shows that
\begin{equation*}
t \le (1-e^{-r})x_1 - \delta(r){x_1\log\log x_1\over\log x_1},
\end{equation*}
which implies that
\begin{equation}
x_1 \ge {t\over1-e^{-r}} + \delta(r){t\log\log t\over\log t}
\label{xorM}
\end{equation}
when $t$ is large enough in terms of $r$. Since $M_t(r)$ equals the
smallest such $x_1$, we see that the right-hand side of the
inequality~(\ref{xorM}) is also a lower bound for $M_t(r)$. This
argument shows that Theorem~\ref{mrtthm} in its entirety follows from
the two propositions.

The reader will have noticed that although Proposition~\ref{estprop}
is stated uniformly for certain rational numbers in the interval $I$,
no use was made of this in deducing Theorems~\ref{denserthm}
and~\ref{mrtthm}. However, we will need the uniformity present in
Proposition~\ref{estprop} in the proofs of Theorems~\ref{finitethm}
and~\ref{lastthm}; it is for this reason that we take the time
to establish the proposition in its current form.

It turns out that the construction used to establish
Proposition~\ref{estprop} proceeds in two stages which, although
similar in spirit, require quite different subsidiary lemmas to
complete. For this reason, we reduce Proposition~\ref{estprop} to the
following two propositions:

\begin{proposition}
Let $I$ be a closed subinterval of $(0,\infty)$. For any real number
$x$ that are sufficiently large in terms of $I$, any rational number
$r\in I$ whose denominator is not divisible by any prime power
exceeding $x\log^{-22}x$, and any integer $R$ satisfying
\begin{equation}
\big| (1-e^{-r})x - \big( 22(1-e^{-r}) - \frac{3r}{e^r} \big)
{x\log\log x\over\log x} - R \big| < \frac r{e^r} {x\log\log
x\over\log x},
\label{Rrange}
\end{equation}
there exists a set\/ \2R of integers satisfying:
\begin{enumerate}
\item[{\rm(i)}]\2R is contained in $[x/2e^r,x]$;
\item[{\rm(ii)}]$|\2R|=R$;
\item[{\rm(iii)}]if $\displaystyle\frac ab = r-\sum_{n\in\2R}\frac1n$
in lowest terms, then $\displaystyle \frac1{\log x} < \frac ab < 1$
and $P^*(b)\le x^{1/5}$.
\end{enumerate}
\label{bigprop}
\end{proposition}
\noindent(As one might think, the constants in the inequality
(\ref{Rrange}), other than the initial $(1-e^{-r})$, are somewhat
arbitrary and chosen simply for convenience during the proof.)

\begin{proposition}
Let $y$ be a sufficiently large real number, and let $a/b$ be a
rational number satisfying $1/\log y<a/b<1$ and $P^*(b)\le y$. Then
there is a set\/ \2S of integers satisfying:
\begin{enumerate}
\item[{\rm(i)}]\2S is contained in $[1,2y^4]$;
\item[{\rm(ii)}]$|\2S|=2\pi^*(y)$;
\item[{\rm(iii)}]$\displaystyle\frac ab = \sum_{n\in\2S} \frac1n$.
\end{enumerate}
\label{smallprop}
\end{proposition}

To see how Propositions \ref{bigprop} and \ref{smallprop} together
imply Proposition~\ref{estprop}, we fix a closed interval
$I\subset(0,\infty)$, an integer $t>T(I)$ where $T(I)$ is a positive
constant that is sufficiently large in terms of $I$, and a rational
number $r\in I$ whose denominator is not divisible by any prime power
exceeding $t\log^{-22}t$. We define
\begin{equation}
x = \frac t{1-e^{-r}} + \big( {22-(22+3r)e^{-r} \over (1-e^{-r})^2}
\big) {t\log\log t\over\log t}
\label{xdef}
\end{equation}
and $R=t-2\pi^*(x^{1/5})$. When these values of $x$ and $R$ are
substituted into the inequality (\ref{Rrange}), the left-hand side has
order of magnitude $t(\log\log t)^2/\log^2t$ after simplification,
while the right-hand side has order of magnitude $t\log\log t/\log t$;
therefore the inequality~(\ref{Rrange}) holds as long as $T(I)$ is
large enough. Certainly $x\log^{-22}x\ge t\log^{-22}t$ as well if
$T(I)$ is large enough. We may therefore apply
Proposition~\ref{bigprop} to obtain a set \2R of integers and a
rational number $a/b$ satisfying properties (i)--(iii) of
Proposition~\ref{bigprop}. With this rational $a/b$, we may then apply
Proposition~\ref{smallprop} with $y=x^{1/5}$ as long as $T(I)$ is
large enough (since $x$ and $y$ are functions of $t$), obtaining a set
\2S satisfying properties (i)--(iii) of that proposition.

We now set $\2E=\2R\cup\2S$. Because $2y^4=2x^{4/5}<x/2e^r$ if $T(I)$
is large enough, it follows from the two properties (i) that \2R and
\2S are disjoint and that the integers in \2E do not exceed
$x$. Moreover, the two properties (ii) imply that $|\2E|=t$, while the
two properties (iii) imply that $r=\sum_{n\in\2E}1/n$. By the
definition~(\ref{xdef}) of $x$, we see that this set \2E satisfies all
the properties required for the conclusion of
Proposition~\ref{estprop}, and so that proposition does indeed follow
from Propositions~\ref{bigprop} and~\ref{smallprop}.

In a sense, we have separated the desire to have an Egyptian fraction
with many terms from the desire to have a specific number of unit
fractions that add to $r$. Proposition \ref{bigprop} yields an
Egyptian fraction with many terms, but one whose sum is only an
approximation to $r$; while Proposition \ref{smallprop} yields an
exact Egyptian fraction representation of $a/b$ with a specified
number of terms, but with that number of terms rather small compared
to the size of the largest denominator.

The broad idea of the proofs of these two propositions uses the
general strategy employed by the author in \cite{Mar:DEF} and by Croot
in \cite{Cro:OSQoEaGaEF}: a collection of unit fractions whose sum is
relatively close to the target rational number is constructed, and
then for each prime power that appears in the denominator of this sum
but not in the denominator of the target rational, a few terms are
omitted or added so that the modified sum is no longer divisible by
that prime power. When all of the unwanted prime powers have been
eradicated in this way, estimates on the number and sizes of the
omitted or added terms are used to show that the resulting sum must
exactly equal the target rational. Proposition \ref{bigprop} is used
to evict the larger unwanted prime powers, while Proposition
\ref{smallprop} is used to evict the smaller prime powers.

To summarize the achievements of this section, we have reduced
Theorems \ref{denserthm} and \ref{mrtthm} to establishing
Propositions~\ref{bestprop},~\ref{bigprop}, and~\ref{smallprop}. These
three propositions will be the subjects of
Sections~\ref{verylargesec},~\ref{crootpropsec},
and~\ref{smallproofsec}, respectively. To complete the outline of the
rest of this paper, we mention that Theorem~\ref{finitethm} will be
established in Section~\ref{anothersec} and that Theorem~\ref{lastthm}
will be established in Section~\ref{yetanothersec}.

\section{The Very Large Prime Powers}\label{verylargesec}\noindent In this section we establish Proposition~\ref{bestprop}, which was
used to show that Theorems~\ref{denserthm} and~\ref{mrtthm} are best
possible. The strategy is to quantify the observation that a tiny
multiple of a very large prime (or prime power) cannot appear in an
Egyptian fraction representation of a given rational number, and then
to calculate the effect that this restriction has on the possible
number of terms in such a representation.

\begin{lemma}
Let $x$ be a positive real number. Suppose that $x_1,\dots,x_t$ are
distinct positive integers not exceeding $x$ and that $p$ is a prime
dividing at least one of the $x_i$. If $p$ does not divide the
denominator of $\sum_{i=1}^t1/x_i$ expressed in lowest terms, then
$p\ll x/\log x$.
\label{bestposslem}
\end{lemma}

\noindent We remark that a slightly modified version of this lemma
could be established for prime powers $q$ rather than merely for
primes $p$, but this formulation suffices for our purposes.
\vskip8pt

\begin{proof}
Because the integers $x_i$ that are not divisible by $p$ do not affect
whether the denominator of $\sum_{i=1}^t1/x_i$ is divisible by $p$, we
may assume that $p$ divides all of the $x_i$. Set $w_i=x_i/p$ for each
$1\le i\le t$, and set $\lambda=\lcm\{w_1,\dots,w_t\}$. Then
\begin{equation}
\sum_{i=1}^t \frac1{x_i} = \bigg( \sum_{i=1}^t {\lambda p\over x_i}
\bigg) \big/ \lambda p,
\label{s12}
\end{equation}
where each summand is an integer by the definition of $\lambda$. We
are assuming that the denominator of $\sum_{i=1}^t1/x_i$ is not
divisible by $p$ when reduced to lowest terms. Consequently, if
\begin{equation}
N = \sum_{i=1}^t {\lambda p\over x_i} = \lambda \sum_{i=1}^t {1\over
w_i}  \label{lameq}
\end{equation}
is the numerator of the fraction on the right-hand side of equation
(\ref{s12}), then $N$ must be a positive multiple of $p$; in
particular, $N\ge p$ and thus $\log N\ge\log p$.

On the other hand, the collection $\{w_i\}$ is a subset of the
integers not exceeding $x/p$. If we define
$L(x)=\lcm\{1,2,\dots,\floor x\}$, so that
\begin{equation}
L(x)=\prod_{p^\nu\le x}p=\exp\bigg( \sum_{p^\nu\le x}\log p \bigg) \le
e^{2x}  \label{cheb}
\end{equation}
by the prime number estimate of Chebyshev, then $\lambda\le L(x/p)$
and so equation (\ref{lameq}) implies
\begin{equation*}
N \le L\big( \frac xp \big) \sum_{w\le x/p} \frac1w \le L\big( \frac
xp \big) \big( \log\frac xp + 1 \big).
\end{equation*}
But then
\begin{equation*}
\log p \le \log N \le \log L\big( \frac xp \big) + O\big(
\log\log\frac xp\big) \ll \frac xp
\end{equation*}
by the estimate (\ref{cheb}), and so $p\ll x/\log x$ as claimed.
\end{proof}

Next we establish an elementary lemma that provides asymptotic
formul\ae\ for the number of integers free of very large prime (or
prime power) factors and for the sum of the reciprocals of such
integers.

\begin{lemma}
Uniformly for $\sqrt x\le y\le x$ and $0<\alpha<1$, we have
\begin{equation*}
\sum \begin{Sb}\alpha x\le n \le x \\ P(n)>y\end{Sb} 1 =
(1-\alpha)x\log{\log x\over\log y} + O\big( \frac x{\log x} \big)
\end{equation*}
and
\begin{equation*}
\sum \begin{Sb}\alpha x\le n \le x \\ P(n)>y\end{Sb} \frac1n = \log
\alpha^{-1}\log{\log x\over\log y} + O\big( \frac1{\alpha\log x} \big).
\end{equation*}
Both formulas remain valid if $P(n)$ is replaced by $P^*(n)$ in the
conditions of summation.
\label{primesums}
\end{lemma}

\begin{proof}
It is a direct consequence of Mertens' formula for $\sum_{p\le x}1/p$
that
\begin{equation}
\sum_{y<p\le x} \frac1p = \log{\log x\over\log y} + O\big( \frac1{\log
y} \big),
\label{mertens}
\end{equation}
and the same asymptotic formula is true if the summation is taken over
prime powers $q$ rather than merely primes $p$. We note that any
integer $n\le x$ such that $P(n)>y$ can be written as $n=mp$ where
$y<p\le x$, and this representation is unique since $y\ge\sqrt x$.
The first assertion of the lemma then follows from equation
(\ref{mertens}) by writing
\begin{equation*}
\sum \begin{Sb}\alpha x\le n \le x \\ P(n)>y\end{Sb} 1 =
\sum_{y<p\le x} \sum_{\alpha x/p \le m \le x/p} 1 = \sum_{y<p\le x}
\big( \frac{(1-\alpha)x}p + O(1) \big)
\end{equation*}
and invoking Chebyshev's estimate $\pi(x)\ll x/\log x$ to bound the
error term. The second assertion follows in a similar manner from
writing
\begin{equation*}
\sum \begin{Sb}\alpha x\le n \le x \\ P(n)>y\end{Sb} \frac1n =
\sum_{y<p\le x} \frac1p \sum_{\alpha x/p \le m \le x/p} \frac1m =
\sum_{y<p\le x} \frac1p \big( \log\alpha^{-1} + O\big( \frac p{\alpha
x} \big) \big).
\end{equation*}
Because the asymptotic formula (\ref{mertens}) is insensitive to the
inclusion of the proper prime powers, these arguments are equally
valid when $P(n)$ is replaced by $P^*(n)$ in the conditions of
summation.
\end{proof}

\begin{pflike}{Proof of Proposition~\ref{bestprop}:}
Suppose that $r$ is a positive rational number, $x$ is a real number
that is sufficiently large in terms of $r$, and \2E is a set of
positive integers not exceeding $x$ such that $\sum_{n\in\2E}1/n=r$.
Let $C>1$ be a large constant and define
\begin{equation*}
\2A = \{n\le x: P(n)\le Cx/\log x\}.
\end{equation*}
We can assume that $x$ is so large that all of the prime divisors of
the denominator of $r$ are less than $Cx/\log x$. Then if $C$ is
chosen large enough, the set \2E must be contained in \2A by Lemma
\ref{bestposslem}.

Choose $0<\alpha<1$ such that, if we set $\2E'=[\alpha x,x]\cap\2A$,
then $|\2E'|=|\2E|$; in other words, $\2E'$ is the subset of \2A with
cardinality $|\2E|$ whose elements are as large as possible. Then
\begin{equation*}
\begin{split}
|\2E| = |\2E'| &= \sum_{\alpha x\le n\le x} 1 - \sum \begin{Sb}\alpha
x\le n\le x \\ P(n)>Cx/\log x\end{Sb} 1 \\
&= (1-\alpha)x+O(1) - (1-\alpha)x\log\big( {\log x\over\log(Cx/\log x)}
\big) + O\big( \frac x{\log x} \big)
\end{split}
\end{equation*}
by Lemma \ref{primesums} with $y=Cx/\log x$. Since
\begin{equation}
\log\big( {\log x\over\log(Cx/\log x)} \big) = \log\big( 1-{\log\log
x\over\log x} + {\log C\over\log x} \big)^{-1} = {\log\log x\over\log
x} + O\big( \frac1{\log x} \big),
\label{lll}
\end{equation}
we see that
\begin{equation}
|\2E| = (1-\alpha)x \big( 1 - {\log\log x\over\log x} \big) + O\big(
\frac x{\log x} \big).
\label{needalpha}
\end{equation}

On the other hand, the elements of $\2E'$ are by definition at least
as big as the elements of \2E, and so
\begin{equation*}
\begin{split}
r = \sum_{n\in\2E}\frac1n \ge \sum_{n\in\2E'}\frac1n &= \sum_{\alpha
x\le n\le x}\frac1n - \sum \begin{Sb}\alpha x\le n\le x \\
P(n)>Cx/\log x\end{Sb}\frac1n \\
&= \log\alpha^{-1} + O\big( \frac1{\alpha x} \big) - \log
\alpha^{-1}\log{\log x\over\log(Cx/\log x)} + O\big( \frac1{\alpha\log
x} \big),
\end{split}
\end{equation*}
again by Lemma \ref{primesums} with $y=Cx/\log x$. Using equation
(\ref{lll}) again, we see that
\begin{equation*}
r \ge \log\alpha^{-1} \big( 1- {\log\log x\over\log x} \big) + O\big(
\frac1{\alpha\log x} \big),
\end{equation*}
which implies that
\begin{equation*}
\alpha \ge e^{-r} \big( 1 - {r\log\log x\over\log x} + O_r\big(
\frac1{\log x} \big) \big).
\end{equation*}

With this lower bound, equation (\ref{needalpha}) becomes
\begin{equation*}
\begin{split}
|\2E| &\le \big( 1- e^{-r} \big( 1 - {r\log\log x\over\log x} +
O_r\big( \frac1{\log x} \big) \big) \big)x \big( 1 - {\log\log
x\over\log x} \big) + O\big( \frac x{\log x} \big) \\
&= (1-e^{-r})x - (1-e^{-r}(1+r)) {x\log\log x\over\log x} +
O_r\big( \frac x{\log x} \big).
\end{split}
\end{equation*}
Since $e^r>1+r$ for $r>0$, we may choose $\delta(r)$ satisfying
$0<\delta(r)<1-e^{-r}(1+r)$, whence
\begin{equation}
|\2E| \le (1-e^{-r})x - \delta(r) {x\log\log x\over\log x}
\label{Ebestp}
\end{equation}
when $x$ is large enough in terms of $r$.\qed
\end{pflike}

\section{The Large Prime Powers}\label{crootpropsec}\noindent In this section we establish Proposition \ref{bigprop}. As mentioned
in the introduction, the methods in this section are in large part
derived from those in Croot \cite{Cro:OSQoEaGaEF}, albeit with some
modifications necessary for the problem at hand. In particular, our
Lemma \ref{crootlem} below is a direct generalization of
\cite[Proposition~2]{Cro:OSQoEaGaEF}.

Let $\|x\|$ denote the distance from $x$ to the nearest integer. If
$a$ and $n$ are coprime, then let $\bar a\mod n$ denote the integer
$b$ with $0<b<n$ and $ab\equiv1\mod n$. (Often we write simply $\bar
a$ when the modulus is clear from the context, e.g., when the term
$\bar a$ appears in the numerator of a fraction whose denominator is
the modulus.) The following lemma demonstrates that under suitable
conditions on a set of integers \2M, the elements of \2M cannot all be
small compared to a modulus $n$ and yet predominantly have
inverses\mod n that are close to 0\mod n, even when scaled by a
nonzero residue $h$.

\begin{lemma}
Let $n$ be a sufficiently large integer, let $k$ be a positive
integer, and let $B$ and $C$ be positive real numbers with $C$
satisfying $200(\log n/\log\log n)^k<C<n$. Suppose that\/ \2M is a set
of positive integers with cardinality greater than $C$, such that each
element $m$ of\/ \2M is less than $B$ and is the product of $k$
distinct primes not dividing $n$. Then for any $0<h<n$, at least $C/2$
elements $m$ of\/ \2M satisfy
\begin{equation*}
\big\| \frac{h\bar m}n \big\| > {C(\log\log n)^k\over200B\log^k
n},
\end{equation*}
where $\bar m$ denotes the inverse of $m\mod n$.
\label{halfbiglem}
\end{lemma}

\begin{proof}
For $m\in\2M$, define $r_m$ to be the integer satisfying $-n/2<r_m\le
n/2$ and $r_m\equiv h\bar m\mod n$; since $n$ does not divide $h$ and
$(m,n)=1$, we see than $r_m$ is nonzero. Also define
\begin{equation*}
s_m = \frac{mr_m-h}n,
\end{equation*}
so that $s_m$ is an integer satisfying $|s_m|<m/2+1$. Suppose that at
least $C/2$ of the $s_m$ satisfied $|s_m|<C(\log\log n)^k/100\log^k
n$. Then, by the pigeonhole principle, there would be a particular
value $s$ with $|s|<C(\log\log n)^k/100\log^k n$ such that $s_m=s$ for
at least
\begin{equation*}
\frac{C/2}{2C(\log\log n)^k/100\log^k n+1} > 20\big( {\log
n\over\log\log n} \big)^k
\end{equation*}
of the elements $m$ of \2M, by the lower bound on $C$. For each such
$m$, we see that
\begin{equation*}
r_m=\frac{ns_m+h}m = \frac{ns+h}m;
\end{equation*}
and since the $r_m$ are nonzero integers, we see that the nonzero
integer $ns+h$ is divisible by at least $20(\log n/\log\log n)^k$
elements of \2M.

On the other hand, it is well-known that the maximal order of the
number of distinct prime divisors of an integer $m$ is asymptotic to
$\log m/\log\log m$, as achieved by those $m$ that are the product of
all the primes up to about $\log m$; thus when $m$ is sufficiently
large, every integer up to $m$ has less than $2\log m/\log\log m$
distinct prime divisors. Since $|ns+h|<n(C(\log\log n)^k/100\log^k
n+1)<n^2$ by the upper bound on $C$, we see that $ns+h$ has at most
$2\log n^2/\log\log n^2<4\log n/\log\log n$ distinct prime factors
when $n$ is sufficiently large, and so the number of divisors of
$ns+h$ that are the product of $k$ distinct primes does not exceed
\begin{equation*}
\frac1{k!} \big( {4\log n\over\log\log n} \big)^k < 20\big( {\log
n\over\log\log n} \big)^k,
\end{equation*}
contradicting the lower bound for the number of elements of \2M that
divide $ns+h$.

This contradiction shows that at most $C/2$ of the $s_m$ satisfy
$|s_m|<C(\log\log n)^k/100\log^k n$, and so at least $C/2$ of the $s_m$
satisfy the reverse inequality. For these elements $m$, we see that
\begin{equation*}
|r_m| = \frac{|ns_m+h|}m > \frac nm \big( {C(\log\log
n)^k\over100\log^k n} - 1 \big) > {Cn(\log\log n)^k\over200B\log^k n}
\end{equation*}
by the upper bounds for $h$ and $m$ and the lower bound for $C$. But
then by the definition of~$r_m$,
\begin{equation*}
\big\| \frac{h\bar m}n \big\| = \big\| \frac{r_m}n
\big\| > {C(\log\log n)^k\over200B\log^k n},
\end{equation*}
which establishes the lemma.
\end{proof}

The next lemma translates the statement of Lemma \ref{halfbiglem} to
an assertion that under suitable conditions on a set \2M, the
inverses\mod n of the elements of \2M must have a subset whose sum is
congruent to any predetermined residue class\mod n.

\begin{lemma}
Let $n$ be a sufficiently large integer, let $k$ be a positive
integer, and let $B$ be a real number satisfying
\begin{equation}
B > {(\log n)^{(k-1)/2}\over(\log\log n)^{k/2}}.  \label{Bbd}
\end{equation}
Suppose that\/ \2M is a set of integers whose cardinality $C$
satisfies
\begin{equation}
C > {200B^{2/3}(\log n)^{(2k+1)/3}\over(\log\log n)^{2k/3}},  \label{Cbd}
\end{equation}
such that each element $m$ of\/ \2M is less than $B$ and is the
product of $k$ distinct primes not dividing $n$. Then for any residue
class $a\mod n$, there is a subset\/ \2K of\/ \2M such that
\begin{equation}
\sum_{m\in\2K} \bar m \equiv a\mod n.  \label{subsetcond}
\end{equation}
\label{crootlem}
\end{lemma}

\begin{proof}
We begin by remarking, in preparation for applying Lemma
\ref{halfbiglem}, that the hypotheses on $B$ and $C$ ensure that
$C>200(\log n/\log\log n)^k$. Also, if $C\ge n$, then the conclusion
of the lemma holds under the weaker assumption that each element of
\2M is coprime to $n$, by the Cauchy--Davenport--Chowla Theorem (see
for instance Vaughan \cite[Lemma 2.14]{Vau:THLM}, and also \cite[Lemma
2]{Mar:DEF}). Therefore we can assume that $C<n$.

Let $e_n(x)$ denote the complex exponential $e^{2\pi ix/n}$ of period
$n$. If we let $N$ be the number of subsets of \2M satisfying the
condition (\ref{subsetcond}), then by the finite Fourier transform
\begin{equation}
\begin{split}
N = \sum_{\2K\subset\2M} \frac1n \sum_{h=0}^{n-1} e_n\bigg( h\bigg(
\sum_{m\in\2K} \bar m - a\bigg) \bigg) &= \frac1n \sum_{h=0}^{n-1}
e_n(-ha) \prod_{m\in\2M} (1+e_n(h\bar m)) \\
& = \frac{2^C}n + \frac1n \sum_{h=1}^{n-1} e_n(-ha) P_h,
\end{split}
\label{N}
\end{equation}
where we have defined $P_h = \prod_{m\in\2M} (1+e_n(h\bar m))$. Using
the identity $|1+e^{it}|^2=2+2\cos t$ and the inequality $1+\cos 2\pi
t \le 2-8\|t\|^2$, we see that
\begin{equation}
|P_h|^2 = \prod_{m\in\2M} \big( 2 + 2\cos\frac{2\pi h\bar m}n \big)
\le 4^C \prod_{m\in\2M} \big( 1 - 4\big\| \frac{h\bar m}n
\big\|^2 \big).
\label{ph2}
\end{equation}

All of the terms in this product are nonnegative and bounded above by
1; and when $n$ is sufficiently large, by Lemma \ref{halfbiglem} at
least $C/2$ of them are bounded above by
\begin{equation*}
1-4\big( {C(\log\log n)^k\over200B\log^k n} \big)^2
\end{equation*}
when $1\le h\le n-1$. Using this fact in the inequality (\ref{ph2})
along with the bound $1-t\le e^{-t}$, we obtain
\begin{equation*}
|P_h|^2 \le 4^C \big( 1-{C^2(\log\log n)^{2k}\over10000B^2\log^{2k}n}
\big)^{C/2} \le 4^C \exp\big( {-C^3(\log\log
n)^{2k}\over20000B^2\log^{2k}n} \big) < 4^C \exp(-2\log n)
\end{equation*}
by the lower bound (\ref{Cbd}) on $C$, and thus $|P_h|<2^C/n$ when
$1\le h\le n-1$.

From this upper bound, we deduce from equation (\ref{N}) that
\begin{equation*}
\big| N-\frac{2^C}n \big| \le \frac1n \sum_{h=1}^{n-1} |P_h| <
{(n-1)2^C\over n^2},
\end{equation*}
which implies that $N>2^C/n^2$. In particular, there do exist subsets
\2K of \2M of the desired type.
\end{proof}

We remark that a stronger inequality for $B$ than (\ref{Bbd}) might be
needed to ensure the existence of a set \2M with the properties
described in the statement of Lemma \ref{crootlem} (for instance, in
the case $k=1$).

The following lemma is the main tool that will be used in our
recursive construction in the proof of Proposition \ref{bigprop}.

\begin{lemma}
Let $0<\xi<1$ be a real number, and let $x$ be a real number that
is sufficiently large in terms of $\xi$. Let $c/d$ be a rational
number and define $q=P^*(d)$, and suppose that $x^{1/5}\le q\le
x\log^{-22}x$. Then there exists a set\/ \2U of integers satisfying:
\begin{enumerate}
\item[{\rm(i)}]\2U is contained in $[\xi x,x]$;
\item[{\rm(ii)}]$|\2U| \le 200(x/q)^{2/3}\log^3x$;
\item[{\rm(iii)}]for each element $n$ of\/ \2U, $P^*(n)=q$;
\item[{\rm(iv)}]if $\displaystyle\frac{c'}{d'} = \frac
cd-\sum_{n\in\2U}\frac1n$ in lowest terms, then $P^*(d')<q$.
\end{enumerate}
\label{biglem}
\end{lemma}

\begin{proof}
We apply Lemma \ref{crootlem} with $n=q$, $k=4$, $B=x/q$, and $a$ the
residue class of $c\overline{(d/q)}\mod q$. Let \2P be the set of all
primes in the interval $((\xi x/q)^{1/4}, (x/q)^{1/4})$ that do not
divide $q$, and let $\2M_0$ be the set of all integers of the form
$p_1p_2p_3p_4$, where the $p_i$ are distinct elements of
\2P. Certainly each element of $\2M_0$ is between $\xi B$ and
$B$. Since the cardinality of \2P is $\gg
(1-\xi)(x/q)^{1/4}\log^{-1}x$ by the prime number theorem, the
cardinality of $\2M_0$ is $\gg_\xi x/(q\log^4x)$, and so we can choose
a subset \2M of cardinality
\begin{equation*}
C=\floor{200(x/q)^{2/3}\log^3x}
\end{equation*}
when $x$ is sufficiently large in terms of $\xi$, by the upper
bound on $q$.

Since these values of $B$ and $C$ do satisfy the hypotheses
(\ref{Bbd}) and (\ref{Cbd}) of Lemma \ref{crootlem}, we can find a
subset \2K of \2M such that $\sum_{m\in\2K} \bar m \equiv a\mod
q$. Now define
\begin{equation*}
\2U = \{qm: m\in\2K\}.
\end{equation*}
We see immediately that properties (i) and (ii) hold for \2U. Since
$q\ge x^{1/5}$, all the elements of \2P (hence all prime-power
divisors of elements of \2K) are less than $q$, and so property (iii)
holds for \2U as well. As for property (iv), it is clear that
$P^*(d')\le q$, since all of the prime powers dividing $d$ or any of
the elements of \2U are at most $q$. On the other hand, $q$ does not
divide $d'$, since if $q=p^\nu$ then
\begin{equation*}
\begin{split}
\frac{c'q}{d'} = \frac {cq}d - \sum_{n\in\2U}\frac qn &= \frac c{d/q}
- \sum_{n\in\2M}\frac1m \\
&\equiv c\overline{(d/q)} - \sum_{n\in\2M}\bar m \equiv 0\mod p,
\end{split}
\end{equation*}
and so $d'$ is divisible by at most $p^{\nu-1}$ after reducing to
lowest terms. Thus $P^*(d')<q$, which establishes the lemma.
\end{proof}

As a last step in the construction used to establish Proposition
\ref{bigprop}, we will be appending a collection of unit fractions
none of whose denominators contain large prime-power divisors. The
following lemma ensures that there are enough such integers in a
suitable range to accommodate this.

\begin{lemma}
Let $0<\eta<1$ be a real number. There exists a constant
$\delta=\delta(\eta)$ such that, for all real numbers $x$ that are
sufficiently large in terms of $\eta$ and for all pairs of real
numbers $\eta\le\alpha,\ep<1$, there are at least $\delta x$ integers
$n$ in the interval $[\alpha x/2,\alpha x]$ satisfying $P^*(n)\le
x^\ep$.
\label{smoothlem}
\end{lemma}

\begin{proof}
We recall that a $y$-smooth integer is one all of whose prime factors
are at most $y$, so that in particular, an integer $n$ is
$x^\ep$-smooth precisely when it satisfies $P(n)\le x^\ep$. It is
well-known (see for instance Hildebrand and Tenenbaum
\cite{HilTen:IwLPF}) that the number of $y$-smooth integers $n\le x$
is
\begin{equation}
x\rho\big( {\log x\over\log y} \big) + O_\eta\big( \frac x{\log x}
\big)  \label{et20}
\end{equation}
uniformly for $x^\eta\le y\le x$; here $\rho(u)$ is the Dickman
function, which is positive, continuously differentiable on
$[1,\infty)$, and satisfies $|\rho'(u)|\le1$ on that interval. In
particular,
\begin{equation}
\rho\big( {\log\alpha x\over\log x^\ep} \big) = \rho\big( \ep^{-1} +
{\log\alpha\over\ep\log x} \big) = \rho(\ep^{-1}) + O_\eta\big(
{1\over\log x} \big)
\label{rhoalpha}
\end{equation}
uniformly for $\eta\le\alpha,\ep<1$. However, $n$ being $x^\ep$-smooth
is a slightly weaker condition than $P^*(n)\le x^\ep$. The asymptotics
for the counting function of $x^\ep$-smooth integers could be shown to
hold for the number of integers $n\le x$ satisfying $P^*(n)\le x^\ep$
as well, but since we only need a weak lower bound for the number of
such integers we argue as follows.

By equation~(\ref{et20}), using equation~(\ref{rhoalpha}) and the
analogous statement with $\alpha$ replaced by $\alpha/2$, we see that
the number of $x^\ep$-smooth integers between $\alpha x/2$ and $\alpha
x$ is $\alpha\rho(\ep^{-1})x/2 + O_\eta(x/\log x)$. On the other hand,
the integers $n$ that are $x^\ep$-smooth but for which $P^*(n)>x^\ep$
are all divisible by at least one prime power $p^\nu>x^\ep$ with $p\le
x^\ep$. Given an integer $k\ge3$, any integer $n$ that is
$x^\ep$-smooth but for which $P^*(n)>x^\ep$ must either be divisible
by the $k$th power of some prime, or else by the square of some prime
exceeding $x^{\ep/k}$. The number of integers $\alpha x/2\le n\le
\alpha x$ that are divisible by the $k$th power of a prime is
$\alpha(1-\zeta(k)^{-1})x/2 + O(x^{1/k})$, while the number of
integers $n\le x$ that are divisible by the square of a prime
exceeding $x^{\ep/k}$ is at most
\begin{equation*}
\sum_{x^{\ep/k}<p\le x^{1/2}} \frac x{p^2}+1 \ll x^{1-\ep/k}.
\end{equation*}
Therefore, if we choose $k$ so large that $\zeta(k)^{-1} >
1-\rho(\eta^{-1})$, then the number of integers $\alpha x/2\le n\le
\alpha x$ such that $P^*(n)\le x^\ep$ is at least
\begin{equation*}
\big( \frac{\alpha\rho(\ep^{-1})x}2 + O_\eta\big( \frac x{\log x}
\big) \big) - \big( \frac{\alpha(1-\zeta(k)^{-1})x}2 + O(x^{1/k})
\big) - O(x^{1-\ep/k}) > \delta x
\end{equation*}
for some constant $\delta=\delta(\eta)$, as long as $x$ is
sufficiently large in terms of $\eta$.
\end{proof}

\begin{pflike}{Proof of Proposition~\ref{bigprop}:}
Let $I=[m,M]$ be a closed subinterval of $(0,\infty)$, and let $x$ be
a real number that is sufficiently large in terms of $I$. Let $r\in I$
be a rational number whose denominator is not divisible by any prime
power exceeding $x\log^{-22}x$. Let
\begin{equation*}
\alpha=e^{-r},\quad \eta=\min\{e^{-M},1/5\},\quad \hbox{and }
\xi=e^{-m},
\end{equation*}
so that $0<\eta\le\alpha\le\xi<1$. Since we are assuming that $x$ is
sufficiently large in terms of $I$, we can assume in particular that
$x$ is sufficiently large in terms of $\xi$ and $\eta$ when appealing
to Lemmas \ref{biglem} and \ref{smoothlem}.  Define $\2A$ to be the
set of all integers $n$ in $[\alpha x,x]$ such that $P^*(n)\le
x\log^{-22}x$, and set
\begin{equation*}
z=\pi^*(x\log^{-22}x) \quad\hbox{and}\quad z'=\pi^*(x^{1/5}).
\end{equation*}
Let $\{q_1,q_2,\dots\}$ denote the sequence of prime powers in
increasing order, and let $p_i$ denote the prime of which $q_i$ is a
power. (By convention we set $p_0=q_0=1$).

Our strategy is to recursively define a sequence $\{a_i/b_i\}$
$(z+1\ge i\ge z')$ of rationals that increase in size as the index $i$
decreases, such that the largest prime-power divisor of each $b_i$ is
less than $q_i$. The first member $a_{z+1}/b_{z+1}$ will be the
difference between our original $r$ and the sum of the reciprocals of
the elements of $\2A$, and each $a_i/b_i$ will be obtained from
the previous $a_{i+1}/b_{i+1}$ by adding several unit fractions whose
denominators belong to $\2A$. The collection of all elements of
$\2A$ not involved in this construction will almost be the set
\2R described in Proposition \ref{bigprop}. This collection will have
slightly fewer than the desired $R$ elements, but we will rectify the
error simply by appending the appropriate number of integers without
large prime-power factors from the interval $[\alpha x/2,\alpha x]$,
and the resulting collection will be our set~\2R.

Define
\begin{equation*}
{a_{z+1}\over b_{z+1}} = r-\sum_{n\in\2A}\frac1n.
\end{equation*}
We recursively define rationals $\{a_z/b_z$, $a_{z-1}/b_{z-1}$, \dots,
$a_{z'}/b_{z'}\}$ and sets $\{\2R_z, \dots, \2R_{z'}\}$ of integers as
follows. Suppose first that $q_i$ divides $b_i$; then we apply Lemma
\ref{biglem} with $c/d=a_{i+1}/b_{i+1}$ and $q=q_i$. The lemma
requires that $q_i=P^*(b_{i+1})$, and since we are supposing that
$q_i$ divides $b_i$ we need only ensure that $P^*(b_{i+1})\le q_i$. In
the case $i=z$, the inequality $P^*(b_{i+1})\le q_i$ is equivalent to
$P^*(b_{z+1})\le x\log^{-22}x$, which is satisfied by the definitions
of $a_{z+1}/b_{z+1}$ and $\2A$ and the hypothesis that the denominator
of $r$ is not divisible by any prime power exceeding
$x\log^{-22}x$. On the other hand, the inequality $P^*(b_{i+1})\le
q_i$ will be satisfied for smaller values of $i$ by the recursive
construction (as we will see in a moment). Let $\2R_i$ be the set \2U
obtained from applying Lemma \ref{biglem}, and let
\begin{equation*}
{a_i\over b_i} = {a_{i+1}\over b_{i+1}} + \sum_{n\in\2R_i} \frac1n
\end{equation*}
in lowest terms, so that by the lemma, $P^*(b_i)<q_i$ and thus
$P^*(b_i)\le q_{i-1}$ (justifying the claim of the previous sentence).

On the other hand, if $q_i$ does not divide $b_i$, then we simply set
$\2R_i=\emptyset$ and $a_i/b_i=a_{i+1}/b_{i+1}$; since
$P^*(b_{i+1})\le q_i$ by the recursive construction and $q_i$ does not
divide $b_{i+1}$, we see that $P^*(b_i)\le q_{i-1}$.

Notice that each $\2R_i$ is a subset of $\2A$, since $\2R_i$ is either
empty or else (by Lemma \ref{biglem}) is contained in $[\xi
x,x]\subset[\alpha x,x]$, and each element $n$ of $\2R_i$ satisfies
$P^*(n)=q_i\le x\log^{-22}x$. Notice also that the various $\2R_i$ are
pairwise disjoint, again since $P^*(n)=q_i$ for $n\in\2R_i$.

Now set $\2R'=\2A\setminus\bigcup_{i=z'}^z\2R_i$, so that
\begin{equation*}
r={a_{z'}\over b_{z'}} + \sum_{n\in\2R'}\frac1n
\end{equation*}
(here we have used the disjointness of the $\2R_i$). We note that the
cardinality of $\2A$ is
\begin{equation}
\begin{split}
|\2A| = \sum \begin{Sb}\alpha x\le n\le x \\ P^*(n)\le
x\log^{-22}x\end{Sb} 1 &= (1-\alpha)x + O(1) - \sum \begin{Sb}\alpha
x\le n \le x \\ P^*(n)>x\log^{-22}x\end{Sb} \frac1n \\
&= (1-\alpha)x - {22(1-\alpha)x\log\log x\over\log x} + O\big( \frac
x{\log x} \big)
\end{split}
\label{Rscard}
\end{equation}
by applying Lemma \ref{primesums} with $y=x\log^{-22}x$. On the other
hand, $\2R'$ is a subset of $\2A$, while the set
$\2A\setminus\2R'$ is simply the union of all the $\2R_i$ and
thus has cardinality
\begin{equation*}
\bigg| \bigcup_{z'\le i\le z} \2R_i \bigg| = \sum_{z'\le i\le z}
|\2R_i| \le \sum_{x^{1/5}\le q\le x\log^{-22}x} 200\big( \frac xq
\big)^{2/3}\log^3x \ll x\log^{-4}x
\end{equation*}
by Mertens' formula (\ref{mertens}) and partial summation. Therefore
the last expression in equation (\ref{Rscard}) represents the
cardinality of $\2R'$ as well.

By the hypothesis (\ref{Rrange}) on $R$, we see that $R$ exceeds the
cardinality of $\2R'$, but by no more than $(4\alpha\log\alpha^{-1})
x\log\log x/\log x$. Let $\2R''$ be any set of $R-|\2R'|$ integers $n$
from the interval $[\alpha x/2,\alpha x]$, each satisfying $P^*(n)\le
x^{1/5}$; we can find such a set by Lemma \ref{smoothlem} with
$\ep=1/5$, since
\begin{equation*}
(4\alpha\log\alpha^{-1}) {x\log\log x\over\log x} \le
{4x\log\log x\over e\log x} < \delta x
\end{equation*}
when $x$ is sufficiently large in terms of
$\delta=\delta(\eta)=\delta(I)$.

Now set $\2R = \2R'\cup\2R''$, so that $\2R$ is contained in the
interval $[\alpha x/2,x]=[x/2e^r,x]$ and the cardinality of $\2R$ is
precisely $R$, and set
\begin{equation*}
\frac ab = r - \sum_{n\in\2R} \frac1n = {a_{z'}\over b_{z'}} -
\sum_{n\in\2R''} \frac1n.
\end{equation*}
Since $P^*(b_{z'})<q_z'\le x^{1/5}$ by the definition of $z'$, and
$P^*(n)\le x^{1/5}$ for all $n\in\2R''$ by the definition of $\2R''$,
we see that $P^*(b)\le x^{1/5}$. It only remains to show that $1/\log
x<a/b<1$ to establish the proposition.

We have
\begin{equation}
\frac ab = r - \sum_{n\in\2R} \frac1n = \log\alpha^{-1} -
\sum_{n\in\2A} \frac1n + \sum_{z'\le i\le z} \sum_{n\in\2R_i}
\frac1n - \sum_{n\in\2R''} \frac1n.
\label{threesums}
\end{equation}
From the definition of $\2A$, and by Lemma \ref{primesums} with
$y=x\log^{-22}x$, the first sum in the last expression of equation
(\ref{threesums}) is
\begin{equation*}
\sum_{\alpha x\le n\le x} \frac1n - \sum \begin{Sb}\alpha x\le n \le x
\\ P^*(n)>x\log^{-22}x\end{Sb} \frac1n = \log\alpha^{-1} + O\big(
\frac1{\alpha x} \big) - {22\log \alpha^{-1}\log\log x\over\log x} +
O\big( \frac1{\alpha\log x} \big).
\end{equation*}
The double sum in equation (\ref{threesums}) is at most
\begin{equation*}
\sum_{x^{1/5}\le q\le x\log^{-22}x} \big( 200\big( \frac xq
\big)^{2/3}\log^3x \big) \frac1{\alpha x} \ll \frac1{\alpha\log^4x},
\end{equation*}
and the last sum in equation (\ref{threesums}) is nonnegative and
at most
\begin{equation*}
(R-|\2R'|) \frac2{\alpha x} \le {(4\alpha\log\alpha^{-1}) x\log\log
x\over\log x} \frac2{\alpha x} = {8\log\alpha^{-1}\log\log x\over\log
x}.
\end{equation*}
Consequently, equation (\ref{threesums}) implies the inequalities
\begin{equation*}
{14\log\alpha^{-1}\log\log x\over\log x} + O\big( \frac1{\eta\log x}
\big) \le \frac ab \le {22\log\alpha^{-1}\log\log x\over\log x} +
O\big( \frac1{\eta\log x} \big)
\end{equation*}
(since $\alpha\ge\eta$), which certainly implies that $1/\log
x<a/b<1$ when $x$ is sufficiently large in terms of $\eta$. This
establishes Proposition~\ref{bigprop}.\qed
\end{pflike}

\section{The Small Prime Powers}\label{smallproofsec}\noindent In this section we establish Proposition~\ref{smallprop}. We are now
concerned more with having precise control over the number of terms in
our Egyptian fractions than with sharply bounding the sizes of their
denominators, as opposed to the case when we considered
Proposition~\ref{bigprop} in the previous section. While the lemmas in
this section appear in Croot~\cite{Cro:OSQoEaGaEF}, we provide proofs
for the sake of completeness and because we state the lemmas in
somewhat different forms.

\begin{lemma}
Let $p^\nu\ge5$ be a power of an odd prime $p$, and let $a$ be any
integer. There exist integers $m_1$ and $m_2$, satisfying
$(p^\nu-3)/2\le m_1<m_2<p^\nu$ and $p\dnd m_1m_2$, such that $\bar
m_1+\bar m_2\equiv a\mod p$.
\label{unetun}
\end{lemma}

\begin{proof}
Assume first that $p\ge5$. Consider the set $\2M=\{p^\nu-(p+3)/2\le
m<p^\nu\}$ of $(p+3)/2$ integers, none of which is a multiple of
$p$. Define two sets of residues\mod p:
\begin{equation*}
\2M_1=\{\bar m:m\in\2M\} \hbox{ and } \2M_2=\{a-\bar m:m\in\2M\},
\end{equation*}
where $\bar m$ denotes the multiplicative inverse of $m\mod p$. Both
$\2M_1$ and $\2M_2$ have $(p+3)/2$ distinct elements\mod p, and each
$\2M_i$ is a subset of the $p$ residue classes\mod p, so by the
pigeonhole principle there must be at least three residue classes $m$
common to $\2M_1$ and~$\2M_2$.

For any such $m$, if we let $m_1=\bar m$ and $m_2=\overline{a-m}$,
then each $m_i$ is in \2M by the definitions of the $\2M_i$, and $\bar
m_1+\bar m_2\equiv a\mod p$. Furthermore, there is precisely one
$m\mod p$, namely $m\equiv\bar2a\mod p$, such that $m_1\equiv m_2\mod
p$ when defined this way. Therefore there is at least one pair
$(m_1,m_2)$ of distinct integers in \2M such that $\bar m_1+\bar
m_2\equiv a\mod{p^\nu}$, and we can assume that $m_1<m_2$ by relabeling
if necessary. Since $p^\nu-(p+3)/2\ge(p^\nu-3)/2$, this establishes the
lemma when $p\ge5$.

On the other hand, if $p=3$ then we must have $\nu\ge2$, and the lemma
can be shown to hold by letting $(m_1,m_2)$ equal $(3^\nu-2,3^\nu-1)$,
$(3^\nu-4,3^\nu-1)$, or $(3^\nu-5,3^\nu-2)$, according to whether $a$ is
congruent to 0, 1, or 2\mod 3.
\end{proof}

The following lemma is one of the two main tools used in our recursive
construction in the proof of Proposition~\ref{smallprop}. This lemma
allows us to control all but the smallest prime powers that can appear
in the denominators of the rational numbers to be constructed.

\begin{lemma}
Let $q\ge4$ be a prime power and let $c/d$ be a rational number with
$P^*(d)\le q$. There exists a set\/ $\2U$ of integers satisfying:
\begin{enumerate}
\item[{\rm(i)}]$\2U$ is contained in $[q^2/5,q^2]$;
\item[{\rm(ii)}]$|\2U|=2$ if $q$ is odd, while $|\2U|=0$ or $1$ if
$q$ is even;
\item[{\rm(iii)}]for each element $n$ of\/ $\2U$, $P^*(n)=q$;
\item[{\rm(iv)}]if $\displaystyle\frac{c'}{d'} = \frac
cd-\sum_{n\in\2U}\frac1n$ in lowest terms, then $P^*(d')<q$.
\end{enumerate}
\label{medlem}
\end{lemma}

\begin{proof}
First assume that $q$ is odd. We apply Lemma \ref{unetun} with $p^\nu=q$
and
\begin{equation*}
a = \begin{cases}c\overline{(d/q)}\mod p, &\rmif q\hbox{ divides }d,
\\ 0\mod p, &\rmif q\hbox{ does not divide }d,\end{cases}
\end{equation*}
finding two distinct integers $m_1$ and $m_2$ in the range
$[(q-3)/2,q)$ such that $\bar m_1+\bar m_2\equiv a\mod p$. Let
$\2U=\{qm_1,qm_2\}$. Then properties (i)--(iii) are easily seen to
hold (the first because $(q-3)/2\ge q/5$ for $q\ge5$), and property
(iv) holds because of the congruence\mod p satisfied by $m_1$ and
$m_2$.

On the other hand, if $q$ is even then $q=2^\nu$ for some $\nu\ge2$.
If $2^\nu$ does not divide $d$ then $P^*(d)<q$ already, and we simply
put $\2U=\emptyset$. If $2^\nu$ does divide $d$, then we easily check
that the set $\2U=\{2^\nu(2^\nu-1)\}$ satisfies properties (i)--(iv).
\end{proof}

The following lemma is the second of the two main tools used in our
recursive construction in the proof of
Proposition~\ref{smallprop}. This lemma allows us to control the
smallest prime powers that can appear in the denominators of the
rational numbers to be constructed. Recall that
$L(x)=\lcm\{1,2,\dots,\floor x\}$, as was defined during the proof of
Lemma \ref{bestposslem}.

\begin{lemma}
Let $c/d$ be a rational number and define $q=P^*(d)$, where $q$ is a
power of the prime $p$. There exists an integer $n$ satisfying:
\begin{enumerate}
\item[{\rm(i)}]$L(q)/(p-1) \le n \le e^{2q}$;
\item[{\rm(ii)}]$P^*(n)=q$;
\item[{\rm(iii)}]if $\displaystyle\frac{c'}{d'} = \frac
cd-\frac1n$ in lowest terms, then $P^*(d')<q$.
\end{enumerate}
\label{smalllem}
\end{lemma}

\begin{proof}
Define $a$ to be the residue class of $c(L(q)/q)\overline{(d/q)}\mod
p$, where $1\le a\le p-1$, and let $n=L(q)/a$. Then properties (i) and
(ii) hold by the definition of $L$ and subsequent remarks, while
property (iii) holds by the choice of $a$.
\end{proof}

\begin{pflike}{Proof of Proposition~\ref{smallprop}:}
Define $z=\pi^*(y)$. Our strategy is very similar to the strategy of
the proof of Proposition~\ref{bigprop}. We recursively define a
sequence $\{a_i/b_i\}$ $(z+1\ge i\ge1)$ of rationals that decrease in
size as the index $i$ decreases, such that the largest prime-power
divisor of each $b_i$ is less than $q_i$. The last member of this
sequence, $a_0/b_0$, will be an integer, and we will show that it must
be zero by bounding its absolute value. The first member
$a_{z+1}/b_{z+1}$ will be our original $a/b$, and each $a_i/b_i$ will
be obtained from the previous $a_{i+1}/b_{i+1}$ by subtracting two
unit fractions (except in a few cases where we subtract only one or
none at all), and the collection of the denominators of all these unit
fractions will almost be the desired set \2S. This collection will
have slightly fewer than the desired $2z$ elements, but we will
rectify the error with a simple modification of the splitting identity
(\ref{split}), and the resulting collection will be our set \2S.

Define $a_{z+1}/b_{z+1}=a/b$ and let $y'=\log y$ and
$z'=\pi^*(y')$. We recursively define rationals $\{a_z/b_z, \dots,
a_1/b_1\}$ and sets $\{\2S_z, \dots, \2S_1\}$ of integers as
follows. If $z'<i\le z$, then we apply Lemma \ref{medlem} with $q=q_i$
and $c/d=a_{i+1}/b_{i+1}$. The requirement of the lemma that
$P^*(b_{i+1})\le q_i$ is satisfied for $i=z$ by the hypothesis of the
proposition and the definition of $z$, and it will be satisfied for
smaller values of $i$ by the recursive construction (as we will see in
a moment). Let $\2S_i$ be the set \2U obtained from applying Lemma
\ref{medlem}, and let
\begin{equation*}
{a_i\over b_i} = {a_{i+1}\over b_{i+1}} - \sum_{n\in\2S_i} \frac1n
\end{equation*}
in lowest terms, so that by the lemma, $P^*(b_i)<q_i$ and thus
$P^*(b_i)\le q_{i-1}$ (justifying the claim of the previous sentence).

If instead $1\le i\le z'$, then we check whether $q_i$ divides
$b_{i+1}$. If not, then we simply set $\2S_i=\emptyset$ and
$a_i/b_i=a_{i+1}/b_{i+1}$; since $P^*(b_{i+1})\le q_i$ by the
recursive construction as before and since $q_i$ does not divide
$b_{i+1}$, we see that $P^*(b_i)\le q_{i-1}$. On the other hand, if
$q_i$ does divide $b_{i+1}$, then we apply Lemma \ref{smalllem} with
$c/d=a_{i+1}/b_{i+1}$ and $q=q_i$. Let $\2S_i$ be the set $\{n\}$
where $n$ is the integer obtained from applying Lemma \ref{smalllem},
and let
\begin{equation*}
{a_i\over b_i} = {a_{i+1}\over b_{i+1}}-\frac1n
\end{equation*}
in lowest terms, so that by the lemma, $P^*(b_i)<q_i$ and thus
$P^*(b_i)\le q_{i-1}$.

Now set $\2S'=\bigcup_{i=1}^z\2S_i$, so that
\begin{equation*}
\frac ab = {a_0\over b_0} + \sum_{n\in\2S'}\frac1n.
\end{equation*}
We claim that in fact $a_0/b_0=0$. It is certainly an integer since
its denominator $b_0$ satisfies $P^*(b_0)\le q_0=1$ by construction,
and it is less than 1 by the hypothesis that $a/b<1$. Moreover, each
element of each $\2S_i$ is at least $q_i^2/5$ if $i>z'$, and at least
$L(q_i)/(p_i-1)$ if $i\le z'$, since these sets resulted from applying
Lemmas \ref{medlem} and \ref{smalllem}, respectively. Therefore, since
$|\2S_i|\le2$ if $i<z'$ and $|\2S_i|\le1$ if $i\le z'$,
\begin{equation*}
\begin{split}
{a_0\over b_0} = \frac ab - \sum_{n\in\2S'} \frac1n &> \frac1{\log y}
- 2\sum_{z'<i\le z} {5\over q_i^2} - \sum_{1\le i\le z'} {p_i-1\over
L(q_i)} \\
&= \frac1{\log y} - \sum_{\log y<q\le y} {10\over q^2} - \sum_{1\le
i\le z'} {p_i-1\over L(q_i)}.
\end{split}
\end{equation*}
Since $L(q_i)/p_i = L(q_{i-1})$, this last sum is a telescoping sum
whose value is $1-1/L(q_{z'})$, as is established by Croot
\cite[Lemma~1]{Cro:OSQoEaGaEF}; in particular, it is less than
1. Also, the penultimate sum is $\ll(\log y\log\log y)^{-1}$ by
Mertens' formula (\ref{mertens}) and partial summation. Therefore
\begin{equation*}
{a_0\over b_0} > \frac1{\log y} - O\big( \frac1{\log y\log\log y}
\big) - 1 > -1
\end{equation*}
when $y$ is sufficiently large. Therefore $a_0/b_0=0$, and
$a/b=\sum_{n\in\2S'}1/n$.

Again by Lemmas \ref{medlem} and \ref{smalllem}, the members of each
$\2S_i$ do not exceed
\begin{equation*}
\begin{cases}
q_i^2, &\rmif z'<i\le z\\
e^{2q_i}, &\rmif 1\le i\le z'
\end{cases}\bigg\}
\le\max\{y^2,e^{2y'}\}=y^2.
\end{equation*}
In addition, if $n$ is a member of $\2S_i$ then $P^*(n)=q_i$, and so
the $\2S_i$ are pairwise disjoint; the cardinality of each $\2S_i$ is
2, except when $1\le i\le z'$ or when $q_i$ is a power of 2, when the
cardinality of $\2S_i$ is 0 or 1. Since there are $\ll\log y$ of these
exceptional values of $i$, we see that $0 \le 2z - |\2S'| \ll \log
y$. Let $n$ be the largest element of $\2S'$ and $m=2z - |\2S'|$, and
define
\begin{equation*}
\2S = (\2S' \setminus \{n\}) \cup \{n+m,n(n+1),(n+1)(n+2),\dots,
(n+m-1)(n+m)\}.
\end{equation*}
Then the cardinality of \2S is exactly $2z=2\pi^*(y)$, and the largest
element of \2S is
\begin{equation*}
(n+m-1)(n+m) \le (y^2+O(\log y))^2 \le 2y^4
\end{equation*}
when $y$ is sufficiently large. Moreover, since the identity
\begin{equation*}
\frac1n = \frac1{n+m} + \frac1{n(n+1)} + \frac1{(n+1)(n+2)} + \dots +
\frac1{(n+m-1)(n+m)}
\end{equation*}
is valid for any positive integers $m$ and $n$, we see also that
$\sum_{n\in\2S} 1/n = \sum_{n\in\2S'} 1/n = a/b$. Therefore \2S
satisfies all of the properties required by Proposition
\ref{smallprop}.
\end{pflike}

\section{The Finiteness of $\2L_j(r)$ for $j\ge2$}\label{anothersec}\noindent In this section we establish Theorem \ref{finitethm}. For the reader's
convenience we recall the definition of the sets $\2L_j(r)$ under
consideration:
\begin{multline*}
\2L_j(r) = \{ x\in{\bf Z},\, x>r^{-1}\colon \hbox{there do not exist }
x_1,\dots,x_t\in{\bf Z},\, x_1>\dots>x_t\ge1 \\
\hbox{ with } { \textstyle\sum_{i=1}^t 1/x_i=r \hbox{ and }
x_j=x } \},
\end{multline*}
so that $\2L_j(r)$ is the set of numbers that {\it cannot\/} be the
$j$th-largest denominator in an Egyptian fraction representation of
$r$. We will make use of the following two lemmas, the first of which
is a simple consequence of Proposition~\ref{estprop} stated in a more
convenient form.

\begin{lemma}
Let $I$ be a closed subinterval of $(0,\infty)$. There exists a
positive real number $X(I)$ such that, for all real numbers $x>X(I)$
and all rational numbers $r=a/b\in I$ for which $P^*(b)<x\log^{-23}x$,
there is a set $\2E$ of positive integers not exceeding $x$ such that
$\sum_{n\in\2E}1/n=r$.
\label{weneedmorelem}
\end{lemma}

\begin{proof}
This follows immediately from Proposition~\ref{estprop} if we set
\begin{equation*}
t = \big\lceil (1-e^{-r})x - {C(I)x\log\log x\over\log x} \big\rceil,
\end{equation*}
where $C(I)$ is a constant that is chosen so large that the right-hand
side of the inequality~(\ref{OIsize}) is less than $x$, and note
that both $t>T(I)$ and $t\log^{-22}t>x\log^{-23}x$ will be true as
long as $X(I)$ is large enough.
\end{proof}

\begin{lemma}
There exists a positive constant $k_0$ such that, for any integer
$k>k_0$, there exists a positive integer $K\equiv-1$\mod k such that
$P^*(K)<k\log^{-24}k$.
\label{allweneedlem}
\end{lemma}

\noindent Because of its length and technical nature, we defer the
proof of Lemma~\ref{allweneedlem} until the end of this
section. Assuming this lemma to be true, we may now proceed with a
proof of Theorem~\ref{finitethm}.

\vskip12pt
\begin{pflike}{Proof of Theorem~\ref{finitethm}:}
We begin by showing that $\2L_2(r)$ is finite for any $r$. Let $r=a/b$
be a positive rational number, and set $I=[r/2,r]$. Let $k$ be any
integer satisfying $k\log^{-24}k>P^*(b)$ and $k>\max\{4/r,X(I),k_0\}$,
where $X(I)$ and $k_0$ are the constants described in
Lemmas~\ref{weneedmorelem} and~\ref{allweneedlem}, respectively. We
claim that there exists an Egyptian fraction representation of $r$
whose second-largest denominator is $k$, and hence that
$k\notin\2L_2(r)$.

To see this, let $K$ be a positive integer such that $K\equiv-1\mod k$
and $P^*(K)<k\log^{-24}k$, as guaranteed by Lemma \ref{allweneedlem}.
Let $x=k/\log k$ and define the rational number $r'=a'/b'$ by
\begin{equation*}
r'=r-{1\over k}-{1\over Kk}.
\end{equation*}
We have that $1/k+1/Kk<2/k<r/2$, and hence $r'$ is in $I$. Also,
\begin{equation*}
\frac1k+{1\over Kk} = {K+1\over Kk} = {(K+1)/k\over K},
\end{equation*}
where the numerator is an integer since $K\equiv-1\mod k$, and so
\begin{equation*}
P^*(b')\le\max\{P^*(b),P^*(K)\} \le k\log^{-24}k < x\log^{-23}x.
\end{equation*}
Therefore we can invoke Lemma~\ref{weneedmorelem} with $r'$ and $x$ to
produce a set \2E of positive integers not exceeding $x$ such that
$\sum_{n\in\2E}1/n=r'$, whence $\2E\cup\{k,Kk\}$ is the set of
denominators for an Egyptian fraction representation of $r$ with
second-largest denominator equal to $k$. Therefore $k$ is not an
element of $\2L_2(r)$, and since this argument holds for all $k$ that
are sufficiently large in terms of $r$ as specified above (note that
the constant $X(I)$ depends only on $r$), we have shown that
$\2L_2(r)$ is finite.

Now that we know that $\2L_2(r)$ is finite, we can establish
Theorem~\ref{finitethm} in its full strength. First, notice that if
there exists an Egyptian fraction representation of $r$ with the
integer $n$ as its $j$th-largest denominator ($j\ge2$), then by
splitting the term with largest denominator using the identity
(\ref{split}), we easily obtain an Egyptian fraction representation of
$r$ with $n$ as its $(j+1)$st-largest denominator. This shows that
$\2L_2(r)\supset\2L_3(r)\supset\2L_4(r)\supset \dotsb$. In particular,
since $\2L_2(r)$ is finite, it follows that all of the $\2L_j(r)$ are
finite for $j\ge3$.

Furthermore, for every element $n$ of $\2L_2(r)$, we can find an
Egyptian fraction representation of $r-1/n$ using only denominators
exceeding $n$ (there are many ways to do this---one could use a greedy
algorithm, for example). If $r-1/n=1/n_1+\dots+1/n_{j-1}$ is such a
representation, then the representation $r=1/n_1+\dots+1/n_{j-1}+1/n$
of $r$ shows that $n\notin\2L_j(r)$. Since we can find such an integer
$j$ for each $n\in\2L_2(r)$, and since the $\2L_j(r)$ form a nested
decreasing sequence of sets each contained in $\2L_2(r)$, we see that
at some point the sets $\2L_j(r)$ will be empty. This completes the
proof of Theorem \ref{finitethm} (modulo the proof of Lemma
\ref{allweneedlem}).\qed
\end{pflike}

We remark that our proof does not show that $\2L_2(r)\subset\2L_1(r)$,
and indeed this is false in general---in fact it is not even the case
that $\2L_3(r)\subset\2L_1(r)$ always. For example, it is easy to see
that if $p$ is a prime, then $p+1$ cannot be the third-largest
denominator in an Egyptian fraction representation of $r=1/p+1/(p+1)$,
and so $p+1$ is an element of $\2L_3(r)$ (hence of $\2L_2(r)$ as well)
but not $\2L_1(r)$. (On the other hand it seems likely that
$\2L_2(1)\subset\2L_1(1)$, for instance, although this doesn't seem
trivial to show.) We can always convert an Egyptian fraction
representation of $r$ whose largest denominator is some integer $n$
into one whose fourth-largest denominator is $n$, by repeatedly
splitting the term with largest denominator other than $n$ and
examining the various ways in which the term $1/n$ could be duplicated
under this process. In this way one can show that
$\2L_4(r)\subset\2L_1(r)$ for every~$r$. It might be interesting to
try to classify the rational numbers $r$ for which
$\2L_2(r)\not\subset\2L_1(r)$.

We now return to the task of establishing Lemma \ref{allweneedlem}.
One possibility would be to cite an existing result on smooth numbers
in arithmetic progressions in which the modulus of the progression was
allowed to exceed the smoothness parameter, such as a theorem of Balog
and Pomerance \cite{BalPom:DSNAP}, and remove those numbers divisible
by large prime powers in an {\it ad hoc\/} manner. We prefer to
provide a self-contained proof of Lemma \ref{allweneedlem}, one that
nevertheless has ideas in common with the method of Balog and
Pomerance, including a reliance on estimates for incomplete
Kloosterman sums.

The following lemma is still much stronger than we need but nearly the
least we could prove to establish Lemma \ref{allweneedlem}. The author
would like to thank Henryk Iwaniec for a helpful conversation
concerning the proof of this lemma.

\begin{lemma}
Let $k$ be a positive integer and $\ep$ and $x$ positive real numbers
such that $k^{6/7+\ep}<x\le k$. The number of ordered pairs $(m,n)$ of
coprime positive integers less than $x$ such that $mn\equiv-1\mod k$
is
\begin{equation}
{6x^2\over\pi^2k} \prod_{p\mid k} \big( {p\over p+1} \big) +
O_\ep(x^{5/6}k^\ep).
\label{klooeq}
\end{equation}
In particular, such ordered pairs $(m,n)$ exist when $k$ is
sufficiently large in terms of $\ep$.
\label{kloolem}
\end{lemma}

\noindent Of course, Lemma \ref{allweneedlem} follows by setting
$\ep=1/14$, say, choosing $x$ such that $k^{13/14}<x<k\log^{-24}k$,
and letting $K=mn$ where $(m,n)$ is one of the pairs whose existence
is ensured by Lemma \ref{kloolem} when $k$ is sufficiently large. The
restriction that $m$ and $n$ be coprime causes a little more trouble
than one might think; without it, one could simplify the proof and
easily obtain an error term of $xk^{-1/4+\ep}$ in the asymptotic
formula (\ref{klooeq}), whence $x$ could be taken as small as
$k^{3/4+\ep}$ in the hypotheses of the lemma. In any case, even a
hypothesis as weak as $x>k^{1-\delta}$ would be ample for our needs,
so we have not gone to great lengths to make the error terms as small
as possible in the proof.

\vskip12pt\begin{pflike}{Proof of Lemma \ref{kloolem}:}
All of the constants implicit in the $\ll$ and $O$-notation in this
proof may depend on $\ep$. We use $d(n)$ to denote the number of
divisors of $n$ and note that $d(n)\ll n^\ep$ for any positive $\ep$.
We also recall that $\|x\|$ denotes the distance from $x$ to the
nearest integer and that $e_k(x)=e^{2\pi ix/k}$.

Let $N$ denote the number of ordered pairs $(m,n)$ of coprime positive
integers less than $x$ such that $mn\equiv-1\mod k$, so that
\begin{equation*}
N = \mathop{\sum_{m\le x}\sum_{n\le x}} \begin{Sb}(m,n)=1 \\
mn\equiv-1\mod k\end{Sb} 1 = \mathop{\sum_{m\le x}\sum_{n\le
x}}_{mn\equiv-1\mod k} \sum_{d\mid(m,n)} \mu(d) = \sum_{d\le x} \mu(d)
\mathop{\sum_{m\le x/d}\sum_{n\le x/d}}_{mnd^2\equiv-1\mod k} 1
\end{equation*}
by changing $m$ and $n$ to $md$ and $nd$, respectively. Let $y$ and
$z$ be real numbers to be specified later subject to $1\le z\le y\le
x$, and write $N=N_1+O(N_2+N_3)$ where
\begin{equation*}
\begin{split}
N_1 &= \sum_{d\le z} \mu(d) \mathop{\sum_{m\le x/d}\sum_{n\le
x/d}}_{mnd^2\equiv-1\mod k} 1; \\
N_2 &= \sum_{z<d\le y} \mathop{\sum_{m\le x/d}\sum_{n\le
x/d}}_{mnd^2\equiv-1\mod k} 1; \\
N_3 &= \sum_{y<d\le x} \mathop{\sum_{m\le x/d}\sum_{n\le
x/d}}_{mnd^2\equiv-1\mod k} 1 \le \sum_{m\le x/y} \sum_{n\le x/y} \sum
\begin{Sb}y<d\le x \\ mnd^2\equiv-1\mod k\end{Sb} 1.
\end{split}
\end{equation*}

We begin by bounding $N_2$ and $N_3$. The estimation of $N_2$ is
trivial: writing $l=mn$ we have
\begin{equation*}
N_2 \le \sum_{z<d\le y} \sum \begin{Sb}l\le x^2/d^2 \\
ld^2\equiv-1\mod k\end{Sb} d(l) \ll (x^2)^{\ep/2} \sum_{z<d\le y}
\big( {x^2\over d^2k} + 1 \big) \ll \big( {x^2\over zk} + y \big)
k^\ep.
\end{equation*}
As for $N_3$, we use Cauchy's inequality to write
\begin{equation*}
N_3^2 \le {x^2\over y^2} \sum_{m\le x/y} \sum_{n\le x/y} \bigg( \sum
\begin{Sb}y<d\le x \\ mnd^2\equiv-1\mod k\end{Sb} 1 \bigg)^2 =
{x^2\over y^2} \sum_{m\le x/y} \sum_{n\le x/y} \mathop{\sum_{y<d_1\le
x}\sum_{y<d_2\le x}} \begin{Sb}mnd_1^2\equiv-1\mod k \\
mnd_2^2\equiv-1\mod k\end{Sb} 1.
\end{equation*}
The congruence conditions imply that $m$ and $n$ must be coprime to
$k$ and thus that $d_1^2\equiv d_2^2\mod k$. Therefore we can weaken
the conditions on the variables $d_1$ and $d_2$ and make them
independent of $m$ and $n$, yielding the upper bound
\begin{equation}
N_3^2 \le {x^2\over y^2} \sum_{m\le x/y} \sum_{n\le x/y}
\mathop{\sum_{y<d_1\le x}\sum_{y<d_2\le x}}_{d_1^2\equiv d_2^2\mod k}
1 \le {x^4\over y^4} \sum \begin{Sb}|l|\le x^2 \\ k\mid l\end{Sb}
\3{1\le d_1,d_2\le x\colon d_1^2-d_2^2=l}
\label{henryk}
\end{equation}
by writing $l=d_1^2-d_2^2$. The term $l=0$ contributes $\floor x$,
while the remaining terms contribute at most $d(|l|)$ each, since
$d_1$ and $d_2$ are determined by $d_1+d_2$ and $d_1-d_2$ which must
be complementary divisors of $l$. Because $|l|\le x^2\le k^2$, we see
that $d(|l|)\ll k^\ep$. Since the number of terms with $|l|\ne0$ in
the final sum in equation~(\ref{henryk}) is $\ll x^2/k\le x$, the sum
is bounded by $O(xk^\ep)$, and so equation~(\ref{henryk}) gives us
\begin{equation*}
N_3 \ll \big( {x^4\over y^4} xk^\ep \big)^{1/2} \le {x^{5/2}k^\ep\over
y^2}.
\end{equation*}
We have thus shown that
\begin{equation}
N = N_1 + O\big( \big( {x^2\over zk} + y + {x^{5/2}\over y^2} \big)
k^\ep \big).
\label{nn1}
\end{equation}

It remains to evaluate $N_1$. Using the additive characters $e_k(x)$
to detect the condition $mnd^2\equiv-1\mod k$, or equivalently
$(dm,k)=1$ and $\bar m+nd^2\equiv0\mod k$, we have
\begin{equation}
\begin{split}
N_1 &= \sum \begin{Sb}d\le z \\ (d,k)=1\end{Sb} \mu(d) \sum
\begin{Sb}m\le x/d \\ (m,k)=1\end{Sb} \sum_{n\le x/d} \frac1k
\sum_{h\mod k} e_k(h(\bar m+nd^2)) \\
&= \frac1k \sum \begin{Sb}d\le z \\ (d,k)=1\end{Sb} \mu(d) \sum_{h\mod
k} \bigg( \sum \begin{Sb}m\le x/d \\ (m,k)=1\end{Sb} e_k(h\bar m)
\bigg) \bigg( \sum_{n\le x/d} e_k(hnd^2) \bigg).
\end{split}
\label{itremains}
\end{equation}
The terms with $h\equiv0\mod k$ contribute
\begin{equation*}
\begin{split}
\frac1k \sum \begin{Sb}d\le z \\ (d,k)=1\end{Sb} \mu(d) \bigg( \sum
\begin{Sb}m\le x/d \\ (m,k)=1\end{Sb} 1 \bigg) \bigg( \sum_{n\le x/d}
1 \bigg) &= \frac1k \sum \begin{Sb}d\le z \\ (d,k)=1\end{Sb} \mu(d)
\big( {x\phi(k)\over dk} + O(k^{\ep/2}) \big) \big( \frac xd+O(1)
\big) \\
&= {x^2\phi(k)\over k^2} \sum \begin{Sb}d\le z \\ (d,k)=1\end{Sb}
{\mu(d)\over d^2} + O\bigg( {x\over k^{1-\ep/2}} \sum
\begin{Sb}d\le z \\ (d,k)=1\end{Sb} \frac1d \bigg) \\
&= {x^2\phi(k)\over k^2} \bigg( {1\over\zeta(2)} \prod_{p\mid k}
\big( 1-{1\over p^2} \big)^{-1} + O\big( \frac1z \big) \bigg) + O\big(
{x\over k^{1-\ep}} \big).
\end{split}
\end{equation*}
Both error terms are $O(x^2/zk^{1-\ep})$, and so equation
(\ref{itremains}) becomes
\begin{equation}
N_1 = {x^2\phi(k)\over\zeta(2)k^2} \prod_{p\mid k} \big( 1-\frac1{p^2}
\big)^{-1} + O\big( {x^2\over zk^{1-\ep}} + \frac1k \sum
\begin{Sb}d\le z \\ (d,k)=1\end{Sb} T(d) \big),
\label{n1eq}
\end{equation}
where we have defined
\begin{equation*}
T(d) = \sum_{h\not\equiv0\mod k} \bigg( \sum \begin{Sb}m\le x/d \\
(m,k)=1\end{Sb} e_k(h\bar m) \bigg) \bigg( \sum_{n\le x/d} e_k(hnd^2)
\bigg).
\end{equation*}

To estimate the $T(d)$, we make use of the elementary bound
\begin{equation*}
\sum_{n\le x/d} e_k(hnd^2) \ll \big\| {hd^2\over k}
\big\|^{-1}
\end{equation*}
when $k$ does not divide $hd^2$, and also the Weil bound for
incomplete Kloosterman sums (see for instance Hooley \cite[Lemma 4
(Section 2.5)]{Hoo:ASMTN}), which gives as a special case
\begin{equation*}
\sum \begin{Sb}m\le x/d \\ (m,k)=1\end{Sb} e_k(h\bar m) \ll_\ep
k^{1/2+\ep/2} (h,k)^{1/2}.
\end{equation*}
It follows that
\begin{equation}
T(d) \ll k^{1/2+\ep/2} \sum_{h\not\equiv0\mod k} (h,k)^{1/2}
\big\| {hd^2\over k} \big\|^{-1}
\label{Tdbd}
\end{equation}
when $d$ is coprime to $k$. Because $(d,k)=1$ we may reindex this sum
by replacing $hd^2$ with $h$, which doesn't affect the quantity
$(h,k)$. When $1\le h\le k/2$, we note that $\|h/k\|^{-1} =
\|(k-h)/k\|^{-1} = k/h$ and $(h,k)=(k-h,k)$, and so
\begin{equation*}
\sum_{h\not\equiv0\mod k} (h,k)^{1/2} \big\| {hd^2\over k}
\big\|^{-1} = 2k \sum_{1\le h\le k/2} {(h,k)^{1/2}\over h}.
\end{equation*}
By writing $h=h'f$ where $f=(h,k)$, this last sum is easily seen to be
$\ll d(k)\log k$, and $T(d)\ll k^{3/2+\ep}$ follows from this bound
and the estimate (\ref{Tdbd}).

Using this estimate for the $T(d)$ in equation (\ref{n1eq}), the
asymptotic formula (\ref{nn1}) for $N$ becomes
\begin{equation}
N = {x^2\phi(k)\over\zeta(2)k^2} \prod_{p\mid k} \big( 1-\frac1{p^2}
\big)^{-1} + O\big( \big( {x^2\over zk} + zk^{1/2} + y + {x^{5/2}\over
y^2} \big) k^\ep \big).
\label{Nlast}
\end{equation}
We optimize this error term by choosing $z=xk^{-3/4}$ and $y=x^{5/6}$,
in which case the error term is $O((xk^{-1/4}+x^{5/6})k^\ep)$, and the
$x^{5/6}k^\ep$ term dominates since $x\le k$. We also note that
$\zeta(2)=\pi^2/6$ and that
\begin{equation*}
{\phi(k)\over k} \prod_{p\mid k} \big( 1-\frac1{p^2} \big)^{-1} =
\prod_{p\mid k} \big( {p\over p+1} \big),
\end{equation*}
and so the main term in (\ref{Nlast}) is the same as the main term in
(\ref{klooeq}). This establishes the lemma.
\end{pflike}

\section{The Order of Magnitude of $L_1(r;x)$}\label{yetanothersec}\noindent In this section we establish Theorem~\ref{lastthm}. We recall that
$L_1(r;x)$ is the number of integers in $\2L_1(r)$ not exceeding $x$,
that is, the number of integers not exceeding $x$ that {\it cannot\/}
be the largest denominator in an Egyptian fraction representation of
$r$.

\vskip12pt\begin{pflike}{Proof of Theorem~\ref{lastthm}:}
First we establish the lower bound in the inequality~(\ref{twosides}).
Let $r$ be a positive rational number and $x>1$ a real number, set
$y=Cx/\log x$ with $C>1$ a large constant, and suppose that $x$ is so
large that all prime divisors of the denominator of $r$ are less than
$y$. By Lemma~\ref{bestposslem}, if $n\le x$ is the largest
denominator in an Egyptian fraction representation of $r$, it must be
true that $n$ is not divisible by any prime larger than $y$ (when $C$
is chosen large enough). In other words, the set $\2L_1(r)$ contains
all integers $n\le x$ such that $P(n)>y$. The number of integers $n\le
x$ such that $P(n)>y$ is asymptotic to $x\log\log x/\log x$ by
Lemma~\ref{primesums}, which establishes the required lower bound.

To establish the corresponding upper bound, we use a method similar to
the one used in the proof of Theorem~\ref{finitethm} to argue that
$\2L_2(r)$ is finite. Let $r=a/b$ be a positive rational number, and
set $I=[r/2,r]$. Let $x$ be a positive real number and set $x'=x/\log
x$ and $y=x'\log^{-23}x'$, so that $y>x\log^{-24}x$.  We suppose that
$x$ is so large that $y>P^*(b)$ and $x'>\max\{2/r,X(I)\}$, where
$X(I)$ is the constant described in Lemma~\ref{weneedmorelem}. Let $k$
be an integer such that $x'<k\le x$ and $P^*(k)<y$, and define the
rational number $r'=a'/b'$ by $r=r-1/k$.

Since $1/k<1/x'<r/2$ we see that $r'\in I$; moreover,
$P^*(b')\le\max\{P^*(b),P^*(k)\}<y$. Therefore we may invoke
Lemma~\ref{weneedmorelem} with $r'$ and $x'$ to produce a set \2E of
positive integers not exceeding $x'$ such that $\sum_{n\in\2E}1/n=r'$,
whence $\2E\cup\{k\}$ is the set of denominators for an Egyptian
fraction representation of $r$ with largest denominator $k$.

From this argument we deduce that for $x$ sufficiently large, the
elements of $\2L_1(r)$ not exceeding $x$ are all contained in
$\{r^{-1}<n\le x'\} \cup \{r^{-1}<n\le x\colon P^*(n)>x\log^{-24}x\}$,
whose cardinality is
\begin{equation*}
\ll x' + x\log{\log x\over\log y} \ll {x\log\log x\over\log x}
\end{equation*}
by Lemma \ref{primesums}. This upper bound holds when $x$ is large
enough in terms of $r$ as described above, but it will hold for all
$x$ by adjusting the implicit constant (depending on $r$) if
necessary. This establishes the upper bound in the inequality
(\ref{twosides}) and hence Theorem~\ref{lastthm}.\qed
\end{pflike}

We remark that we have actually established the inequalities
\begin{equation*}
{(1+o_r(1)) x\log\log x\over\log x} \le L_1(r;x) \le {(24+o_r(1))
x\log\log x\over\log x}.
\end{equation*}
We did so by showing essentially that an integer $n$ is in $\2L_1(r)$
if and only if $n$ can be written as $n=p^\nu m$ with $m$ less than a
power of $\log p$, though we were unable to pinpoint this power other
than to show that it lies between $1+o(1)$ and $24+o(1)$. With much
more care we could improve the constant 24 to 3 but no further at
present. Nevertheless we speculate that the correct power is 1, i.e.,
that for any fixed $r$ and $\ep>0$ there are only finitely many
integers $n\in\2L_1(r)$ such that, if $P^*(n)=p^\nu$ and $n=p^\nu m$,
then $m\ge\log^{1+\ep}p$. This would imply that the counting function
$L_1(r;x)$ of those integers that cannot be the largest denominator in
an Egyptian fraction representation of $r$ is asymptotic to $x\log\log
x/\log x$ with leading coefficient 1.

\bibliography{denser}
\bibliographystyle{amsplain}
\end{document}